%% file: main.tex
\documentclass[onefignum,onetabnum]{siamart220329}
\usepackage{nccmath}
\usepackage{mathtools}
\usepackage[font=footnotesize]{caption}
\usepackage{subcaption}
\usepackage{multirow}
\usepackage{microtype}
\usepackage{adjustbox}
\usepackage{amsmath, amssymb}
\usepackage{float}
\usepackage{cleveref}
\crefname{algocf}{alg.}{algs.}
\Crefname{algocf}{Algorithm}{Algorithms}
\let\vector\undefined
\newcommand\vector[1]  {\boldsymbol{#1}}
\newcommand\vx{{\vector{x}}}
\newcommand\vtheta{{\vector{\theta}}}

\newcommand\vn{{\vector{n}}}
\newcommand\vk{{\vector{k}}}
\newcommand\vN{{\vector{N}}}

\newcommand\vf{{\vector{f}}}

\newcommand\vb{{\vector{b}}}
\newcommand\vc{{\vector{c}}}
\newcommand\vNbar{{\bar{\vector{N}}}}
\newcommand\vm{{\vector{m}}}
\newcommand\vp{{\vector{p}}}
\newcommand\vzero{{\vector{0}}}

\newcommand\PDd{{\mathbb{P}^D_d}}
\newcommand\PDdone{{\mathbb{P}^D_{d,1}}}
\newcommand\PDdtwo{{\mathbb{P}^D_{d,2}}}

\input{shared}

\ifpdf
\hypersetup{
  pdftitle={A Fast Chebyshev Transform for High-Dimensional Smooth Function Approximation},
  pdfauthor={Dalton Jones, Pierre-David Letourneau, Matthew J. Morse and M. Harper Langston}
}
\fi

\begin{document}

\maketitle

% REQUIRED
\begin{abstract}
     We present the Fast Chebyshev Transform (FCT), a fast, randomized algorithm to compute a Chebyshev approximation of functions in high-dimensions from the knowledge of the location of its nonzero Chebyshev coefficients.
     Rather than sampling a full-resolution Chebyshev grid in each dimension, we randomly sample several grids with varied resolutions and solve a least-squares problem in coefficient space in order to compute a polynomial approximating the function of interest across all grids simultaneously.  We theoretically and empirically show that the FCT exhibits quasi-linear scaling and high numerical accuracy on challenging and complex high-dimensional problems.
     We demonstrate the effectiveness of our approach compared to alternative Chebyshev approximation schemes.  In particular, we highlight our algorithm's effectiveness in high dimensions, demonstrating significant speedups over commonly-used alternative techniques.

    %This paper introduces an algorithm for computing the Chebyshev coefficients of a polynomial in quasi-linear time given knowledge of the locations of its non-zero coefficients. We show how our algorithm can be used to perform Chebyshev interpolation in high dimension efficiently and compare the performance of our method with that of other commonly-used techniques based on least-squares and the Discrete Cosine Transform (DCT). Our results demonstrate significant improvements, especially in high dimension.
\end{abstract}

% REQUIRED
\begin{keywords}
  Chebyshev polynomial, high dimensional approximation, Discrete Cosine Transform, randomized sampling, least squares, fast algorithms
\end{keywords}

% REQUIRED
\begin{MSCcodes}

\end{MSCcodes}

\section{Introduction}

Polynomial interpolation is a fundamental building block of numerical analysis and scientific computing, often serving as the connection between continuous application domains and discrete computational domains.
Further, Chebyshev polynomials are a popular basis for interpolation due to their favorable convergence properties for smooth and analytic functions, their numerical stability, their near-best approximation power \cite{li2004near}, and a simple FFT-based interpolation procedure \cite{trefethen2019approximation,trefethen2000spectral,powell1967maximum}.
As a result, Chebyshev polynomial interpolation has been used in a variety of domains, such as solving partial differential equations \cite{mason1967chebyshev,trefethen2000spectral,deville1985chebyshev}, signal processing \cite{shuman2011chebyshev,kabal1986computation} and optimization \cite{cheng1999nonlinear, sherali2001global}, and they serve as the foundation for the popular numerical computing toolkit \texttt{chebfun} \cite{battles2004extension}.

The standard approach for interpolating a function in $D$ dimensions with a $d$ degree polynomial is to sample the target function on a tensor product grid of $(d+1)^D$ points and solve a linear system for the coefficients of a polynomial $p(\vx)$; i.e., 

\begin{equation}
p(\vx) = \sum_{\vn=\vzero}^{\vector{D}} c_{\vn} T_{\vn}(\vx),
\label{eq:chebyshev_expansion}
\end{equation}
where $\vn = (n_1, \hdots, n_D)$ is a multi-index, $\vx = (x_1, \hdots, x_D) \in \mathbb{R}^D$,\newline $T_\vn(\vx) = \prod_{i=1}^D T_{n_i}(x_i)$, and $T_n(x)$ is the $n$th Chebyshev polynomial. 
However, as the dimension grows, the number of required samples and unknown coefficients grow exponentially, and the problem quickly becomes intractable.

Due to the importance of high-dimensional function approximation, many alternative methods have been developed to overcome this {\em curse of dimensionality}.
A more thorough review of the literature can be found in \cite{trefethen2017cubature}, but we briefly summarize various approaches here.
Low-rank approximations \cite{bebendorf2011adaptive,hashemi2017chebfun,townsend2013extension, oseledets2011tensor,grasedyck2013literature} try to reduce the problem to computing coefficients of a reduced set of separated functions that still approximate the target function well. 
The accuracy of low-rank approximations can be improved by increasing the number of terms in the approximation or by applying the method hierarchically on low-rank submatrices  of the resulting system matrix \cite{bebendorf2007hierarchical,hackbusch2012tensor}.

Sparse grids \cite{bungartz2004sparse, griebel2014fast,griebel2021generalized} reduce the degrees of freedom by using a reduced tensor-product basis formed from a progressively finer hierarchy of functions. 
After removing basis elements with mixed derivative terms below a threshold tolerance, the remaining coefficients decay rapidly, allowing for scaling independent of dimension (up to logarithmic factors).
Adaptive methods have been proposed \cite{conrad2013adaptive}, leveraging the pseudospectral method combined with sparse grids and Smolyak's algorithm \cite{smolyak1963quadrature} to reach even higher dimensions.

Recently, deep neural networks have demonstrated their effectiveness at approximating and interpolating high-dimensional data \cite{devore2021neural}. 
This has become more relevant in light of the ``double-descent'' phenomenon \cite{nakkiran2021deep,belkin2019reconciling}: test accuracy actually continues to improve even after a network overfits the training data.
Some progress has been made to determine the theoretical underpinning of this behavior \cite{allen2019convergence,xie2022overparameterization,belkin2019does}.
Although designed for Gaussian process regression and limited to low dimensions, the most recent work that resembles our proposed approach is \cite{greengard2022equispaced}, though our appraoch is more broadly scoped and applied. 

Standard tensor product representations for Chebyshev interpolation in high-dimensions have largely fallen by the wayside due to their exponential complexity. 
The most relevant improvement to the baseline algorithm, while still leveraging tensor product grids, is a faster DCT algorithm that reduces the constant runtime prefactor \cite{chen2003fast}.
Recently, however, \cite{trefethen2017cubature,trefethen2017multivariate} have shown that Chebyshev coefficients decay rapidly when approximating smooth functions, but anisotropically; that is, the normal tensor product grid underresolves the diagonal of the hypercube$[-1,1]^D$ by a factor of $\sqrt{D}$.
Interpolating with a bounded Euclidean degree instead of bounded total degree will address this anisotropy and recover the usual convergence behavior.

Another important observation from \cite{trefethen2017cubature} is that, regardless of choosing a bounded total or Euclidean degree interpolant, the Chebyshev coefficients above machine precision are contained in a sphere of radius $d$ under the $s$-norm ($s=1$ for bounded total degree, $s=2$ for Euclidean).
An obvious consequence of this is that by increasing the radius of this sphere in coefficient-space, we increase the accuracy of the resulting Chebyshev interpolant. 
In other words, to achieve a given target accuracy, $\epsilon$, we can truncate coefficients outside of this coefficient sphere of radius $d_\epsilon$ without losing accuracy. 
This is algorithmically important because there are exponentially fewer coefficients inside a sphere under the 1-or 2-norm than under the $\infty$-norm (i.e. bounded maximal degree) \cite{smith1989small,trefethen2017multivariate}. 

In this paper, we describe an algorithm, the Fast Chebyshev Transform (FCT), whose runtime is proportional to the number of coefficients in this coefficient-space $s$-norm sphere rather than the size of a $D$-dimensional tensor-product grid.
To do this, we estimate the number of significant coefficients $N$ and sample $L$ distinct tensor product Chebyshev grids with random resolutions in each dimension, such that the product of the resolutions roughly equals $N$.
We then sample the target function at all of these grid points, perform separate $D$-dimensional DCTs on each grid, and solve a least-squares problem in coefficient space with the conjugate gradient (CG) method.
Since the system matrix is well-conditioned, CG converges rapidly.
More importantly, we can apply each matrix-vector multiplication in $O(N)$ time due to the sparsity pattern induced by an aliasing identity derived from Chebyshev polynomials.
The FCT algorithm has an overall complexity of $O(N\log N)$.
We limit our scope to Chebyshev approximation in this work, but FCT can be applied directly to high dimensional cosine approximation of smooth functions with minimal modifications.

The remainder of this paper is structured as follows: In \Cref{sec:preliminary}, we discuss relevant notation and background information.
In \Cref{sec:main}, we discuss the mathematical formulation that underpins the Fast Chebyshev Transform and present the full algorithm.
\Cref{sec:algorithm} discusses the algorithmic complexity and numerical behavior of FCT. Numerical results and comparisons are presented in \Cref{sec:numerical_results}. 
We conclude with \Cref{sec:conclusions}. 

\section{Background and notation}
\label{sec:preliminary}
We aim to approximate a given smooth function $f: \mathbb{R}^D \to \mathbb{R}$ by a tensor product Chebyshev polynomial in the form of \Cref{eq:chebyshev_expansion}.
Bold script will represent vector quantities, such as the real vector $\vx=(x_1, \hdots, x_D)$ and integer multi-indices $\vn = (n_1,\hdots, n_D)$. 
The $n$th Chebyshev polynomial is given by 
\begin{equation}
 T_{n}(x) = \cos( n \, \arccos(x)).
\end{equation}
Chebyshev polynomials are orthonormal with respect to measure $\frac{1}{\sqrt{1-x^2}}, \, x\in[-1,1]$:
\begin{equation}
    \int_{-1}^1 \frac{T_m (x) \, T_n(x) }{\sqrt{1-x^2}} dx=
\left\{
	\begin{array}{ll}
		0  & \mbox{if } n \not = m \\
		1 & \mbox{if } n = m
	\end{array}
\right.
\end{equation}
and obey the recurrence relation 
\begin{equation}
    T_{n+1} = 2xT_n(x) - T_{n-1}(x).
\end{equation}
A multi-dimensional Chebyshev polynomial defined on $\mathbb{R}^D$ is a tensor product of one dimensional Chebyshev polynomials:
\begin{equation}
T_\vn(\vx) = \prod_{i=0}^D T_{n_i}(x_i).
\end{equation}

When approximating a function in $D$ dimensions by polynomials, we often consider the class of polynomials with bounded degree $d$ with respect to some norm as approximation candidates; i.e. members of the set

\begin{equation}
\mathbb{P}^D_{d,s} = \left\{p(\vx) = \sum_\vn c_\vn T_\vn(\vx) \in \PDd \biggm| ||\vn||_{s} \leq d \right\},
\end{equation}
where $\PDd$ is the space of polynomials $p:\mathbb{R}^D \to \mathbb{R}$ with maximum degree $d$, and $||\vn||_{s}$ is the $s$-norm of the vector $\vn$.
We will refer to $\PDdone$, $\PDdtwo$ and $\PDd$, as polynomials with bounded total degree, Euclidean degree and maximum degree respectively. 

Increasing the degree of a polynomial approximation of $f$ increases the approximation accuracy \cite{trefethen2019approximation}. 
This means that to approximate $f$ by a polynomial $p$ with accuracy $\epsilon$, there is a degree $d = d_{\epsilon,s}$ such that $\| f-p \| \leq \epsilon$ if $p$ is in $\mathbb{P}^D_{d,s}$.
Choosing $d$ determines the total number of coefficients of basis elements with bounded exponents in the $s$-norm, i.e. $\|\vn \|_s \leq d$, which define $p$.
We denote this set of indices by
\begin{equation}
\mathcal{S}^s_{d,D} =\{\vn \in \mathbb{N}^D \,\mid\, \|\vn\|_s \leq d \}.
\end{equation}
Known results tell us that
\begin{align}
|\mathcal{S}^1_{d,D}     | &= \binom{d + D}{d} \quad \cite{moore2018monomial}\label{eq:N_1norm}\\
|\mathcal{S}^2_{d,D}     | &= \frac{(\sqrt{\pi}d/2)^D}{\Gamma\left(\frac{D}{2} + 1 \right)}, \quad D \to \infty\quad \cite{smith1989small}\label{eq:N_2norm}\\
|\mathcal{S}^\infty_{d,D}| &= d^D \label{eq:N_infnorm}.
\end{align}
The important aspect of \Cref{eq:N_1norm,eq:N_2norm} is that they grow much slower than \Cref{eq:N_infnorm} as $D$ increases.
This makes them more efficient to work with directly than \Cref{eq:N_infnorm} while still approximating $f$ well.
These values are used to discuss the computational complexity of this work; we typically choose $N = |\mathcal{S}^1_{d,D}|$ or $|\mathcal{S}^2_{d,D}|$. 

\paragraph{Problem setup}
One way to approximate a smooth function $f$ is to sample it on a $D$-dimensional tensor product of 1D Chebyshev grids:
\begin{equation}
    x_{k}= \cos(\theta_{k}),\quad \theta_{k} = \frac{\left(k + \frac{1}{2}\right)\pi}{N+1}, \quad k=0,\hdots, N.
    \label{eq:chebyshev_grid}
\end{equation}
Then, apply a DCT in each dimension.  
This produces a polynomial interpolant of the function samples; since Chebyshev polynomials obey the following relation after a change of variables from $[-1,1]$ to $[0,\pi]$
\begin{equation}
    \label{eq:chebeq}
    T_{n} (x_k) = T_{n} (\cos(\theta_k)) = \cos(n \theta_k),
\end{equation}
\Cref{eq:chebyshev_expansion} becomes a multi-dimensional cosine series.
However, as mentioned above, the size of the tensor-product grid grows exponentially with $D$, limiting this approach's practicality to low dimensions.

Another common approach is to approximate $f$ with a least-squares solution: find $p$ that minimizes $\|f - p\|_2^2$.
The standard way to solve this problem is sample $f$ at a set of points, then form a matrix $B$ with entries
\begin{equation}
    B_{i,j} = T_{\vn^{(j)}} \left(\vx^{(i)}\right),
    \label{eq:least_squares_matrix}
\end{equation}
where $\vn^{(j)}$ and $\vx^{(i)}$ imply an ordering on the multi-indices and sample points. 
After constructing $B$, we solve the minimization problem,
\begin{equation}
   \min_\vc \|B\vc - \vf\|_2^2,
    \label{eq:least_squares_opt}
\end{equation}
where $\vc_i = c_{\vn^{(i)}}$ and $\vf_i = f(\vx^{(i)})$.
This can solved many ways, either with a direct approach, such as solving the normal equations, or by using the singular value decomposition (SVD), or with an iterative solver like the conjugate gradient method (CG) \cite{kershaw1978incomplete}.

The advantage of a least-squares-based approach is that the complexity is in terms of the number of unknown coefficients $N$ of \Cref{eq:N_1norm,eq:N_2norm}, rather than an exponentially growing grid size.
A disadvantage of this approach is the implicit ambiguity in choosing sample point locations:
poorly sampling $f$ can produce a bad approximation or slow convergence, while provable convergence often requires exponentially many more samples \cite{demanet2019stable}.
Another drawback is the relatively large complexity; 
direct methods require $O(N^3)$ work, while the CG method still requires $O(N^2)$ dense matrix-vector multiplications.

\section{Algorithm}
\label{sec:main}
\begin{figure}[!htp]
     \centering
     \includegraphics[width=\textwidth]{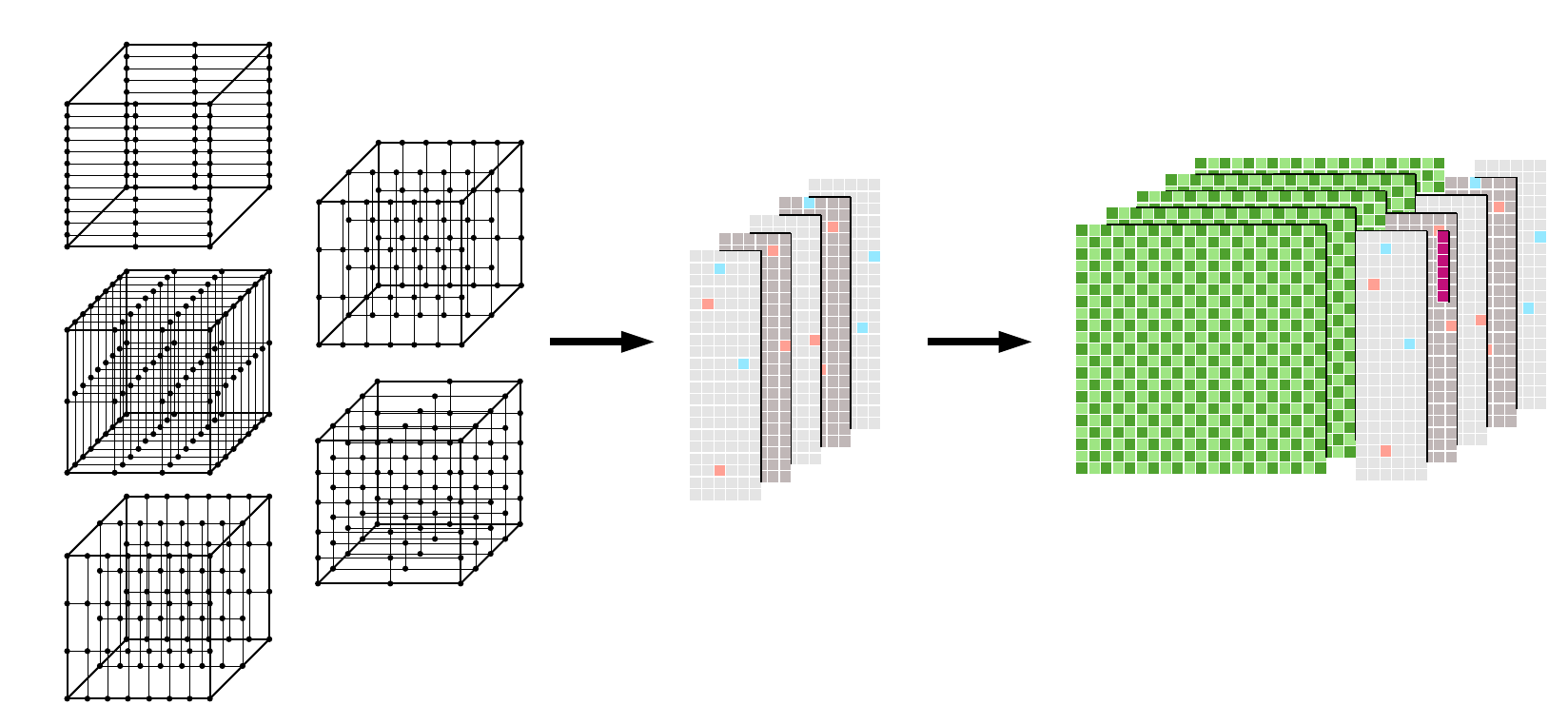}
     \hfill
        \caption{Schematic of the matrix-vector multiplication step in the sparse Fast Chebyshev Transform. {\em Left}:
        First, we sample the high-dimensional input function on a collection of $L$ grids with varied, randomly chosen discretizations in each dimension, the product of each of which roughly equals the estimated number of non-zero coefficients $N$.
        {\em Middle}: Next, for each sampled grid, we compute the non-zero entries (red and blue entries) of the linear system resulting from an interpolation problem on the grid. This is sparse due to aliasing effects of Chebyshev polynomials (see \Cref{lem:aliasing} and  \cite[Chapter 4]{trefethen2019approximation}).
        {\em Right}: Finally, we multiply each sparse matrix with the coefficient vector (magenta) and apply a Discrete Cosine Transform (green matrix) to the result. To solve for unknown polynomial coefficients, we solve a sparse least-squares using an iterative method such as the Calgorithm.}
     \label{fig:fct_schematic}
\end{figure}

To address these shortcomings, we propose an alternative approach: the Fast Chebyshev Transform (FCT). 
The core of the FCT algorithm is a consistent linear system of the form,
\begin{align}
\label{eq:main_linear_system}
\begin{bmatrix}
F^{(1)} & 0 & \hdots & 0 \\
0 & F^{(2)} & \hdots & 0 \\
\vdots & \vdots & \ddots & \vdots \\
 0 & 0 & \hdots & F^{(L)}
\end{bmatrix} 
&\,
\begin{bmatrix}
    \vf^{(1)}\\
    \vf^{(2)}\\
   \vdots \\
   \vf^{(L)}
\end{bmatrix} = 
\begin{bmatrix}
A^{(1)} \\
A^{(2)} \\
\vdots \\
A^{(L)} \\
\end{bmatrix} \, 
\vc
\\
\mathcal{F}\vf &= A\vc.
\end{align}
Each matrix $A^{(\ell)}$ is a sparse  $O( N ) \times N$ matrix, which are concatenated into a large $O(NL)\times N$ matrix with $L$ sparse blocks.
The block matrices $F^{(\ell)}$ are $O(N) \times O(N)$ matrices corresponding to DCTs, $\vc =  [c_{\vn^{(1)}},\hdots,c_{\vn^{(N)}}]^T$ is the vector containing non-zero coefficients of $f$'s polynomial approximation, and $[\vf^{(1)},\hdots, \vf^{(L)} ]^T$ is a vector of $f$ at sample points $\vx^{(i)}$.
To determine the sampling pattern of $f$ and the sparsity pattern of $A$, we construct a collection of $L$ distinct tensor product grids, each containing randomly selected resolutions in each dimension, but whose total number of points is roughly $N$.
We call this collection of grids an \textit{$L$-grid}; an example of one possible $L$-grid for $N=90$ in three dimensions is shown in \Cref{fig:fct_schematic}-Left.

Due to the aliasing properties of Chebyshev polynomials (\Cref{lem:aliasing},  \cite[Chapter 4]{trefethen2019approximation}), the matrices $A^{(l)}$ are sparse and contain $O(N)$ non-zero entries.
These properties, along with fast DCT algorithms \cite{chen1977fast,ahmed1974discrete}, allow us to compute both sides of \Cref{eq:main_linear_system} with $O(LN\log(N))$ time complexity. 
When $L$ is well-chosen, the matrix $A$ has favorable conditioning properties, allowing an iterative solver such as the CG method \cite{kershaw1978incomplete}, to converge rapidly to the least-squares solution.
A schematic of the algorithm is shown in \Cref{fig:fct_schematic}.

In the remainder of this section, we first highlight the aliasing properties of Chebyshev polynomials on $L$-grids and how this leads to the sparsity of $A$ in \Cref{sec:aliasing}. 
In \Cref{sec:sampling_schedule}, we describe how to appropriately choose an $L$-grid.
%Then, in \cref{sec:form_linear_system}, we outline a fast algorithm to construct $A$ from a sequence of sparse Kronecker products. 
We combine each of these pieces in \Cref{sec:algorithm} to describe the full FCT algorithm.

\subsection{Aliasing}
\label{sec:aliasing}
 Chebyshev polynomials evaluated on a Chebyshev grid exhibit a special property known as \textit{aliasing}: a higher degree Chebyshev polynomial realizes the same values as a single lower degree Chebyshev polynomial on a fixed grid size \cite[Chapter 4]{trefethen2019approximation}. 
 When approximating functions with Chebyshev polynomials, this aliasing has another side effect: if a low degree Chebyshev polynomial is zero on a given Chebyshev grid, certain higher degree Chebyshev polynomials will also be zero.
 We leverage this to produce a sparse system matrix in \Cref{eq:main_linear_system}.

 \paragraph{Aliasing in 1D} 
To highlight this relationship, consider the  one-dimensional problem of interpolating a smooth function $f$ by a degree $N$ Chebyshev polynomial on $[-1,1]$.
We know that $f$ can be expressed as a unique Chebyshev series on $[-1,1]$:
\begin{equation}
    f(x) = \sum_{n=0}^\infty c_nT_n(x).
    \label{eq:f_inf_cheb_expansion}
\end{equation}
To approximate this, we want to find a set of $M$ coefficients $c_i$ defining a polynomial $p(x)$ such that 
\begin{equation}
 f(x) \approx p(x) = \sum_{n=0}^M c_nT_n(x), \quad x \in [-1,1], \quad c_n = \frac{2 - \delta_0(n)}{\pi}\int_{-1}^1 \frac{T_n(x)f(x)}{\sqrt{1-x^2}} dx,
\end{equation}
where $\delta_0(n)$ is the dirac delta function centered at $0$. Specifically, we would like to approximate $f$ by a truncated Chebyshev series.
After applying a change of variables $x = \cos(\theta)$ from $[-1,1]$ into $[0,\pi]$, the interpolation problem becomes 
\begin{align}
f(\cos(\theta)) \approx p(\cos(\theta)) &= \sum_{n=0}^M c_n\cos(n\theta), \quad \theta \in [0,\pi], \quad \quad \\
c_n &= \int_{0}^\pi \cos(n\theta)f(\cos(\theta)) d\theta.
\label{eq:poly_approx}
\end{align}

Substituting the polynomial approximation of $f$ into the expression for the $n$th Chebyshev coefficient in \Cref{eq:poly_approx}, we obtain
\begin{align}
c_n =& \int_{0}^\pi \cos(n\theta) \left(\sum_{m=0}^M c_m\cos(m\theta)\right) d\theta \label{eq:cn_equation} \\=& \sum_{m=0}^M c_m\int_{0}^\pi \cos(n\theta) \cos(m\theta) d\theta \label{eq:inner_product_1D}
\end{align}
for $n\leq m$, thanks to the orthogonality of cosines over $[0,\pi]$. If we discretize the integral in \Cref{eq:inner_product_1D} using $N+1$ equispaced points in $[0, \pi]$ (corresponding to a first-kind Chebyshev grid of $N+1$ points in \Cref{eq:chebyshev_grid}), the quadrature becomes
\begin{equation}
b_{n,N} :=\sum_{m=0}^M c_m \left ( \frac{1}{N+1} \, \sum_{k=0}^N \cos(n\theta_k) \cos(m\theta_k) \right ) =   \frac{1}{N+1} \, \sum_{k=0}^N \cos(n\theta_k) \, p \left( \cos(m\theta_k) \right ). 
\label{eq:quad_1D} 
\end{equation}
%where the $\frac{1}{N+1}$ factor is the quadrature weight from the trapezoidal rule. 
In particular, note that if $N \geq M$, the trapezoidal rule is exact for the $(M+1)$-term cosine series in \Cref{eq:inner_product_1D} and $b_{n,N}  = c_n$. However, this is not the case when $N \leq M$ due to \emph{aliasing}. Such aliasing is quantified by the following lemma:

% At first glance, it appears that each coefficient requires quadratic work to compute the sum in \Cref{eq:inner_product_1D}.
% Fortunately, thanks to the orthogonality and aliasing properties of cosine series on equispaced grids, we can quickly estimate the few values of $n$ and $m$ where the quadrature sum is non-zero and approximate the sum analytically. 
% We summarize the main result in the following lemma and provide a proof in 
% \cref{sec:lemma_proof}.
\begin{lemma}
\label{lem:aliasing}
 For every positive integer $N$,  
 \begin{align}
    \frac{1}{N+1} \, &\sum_{k=0}^{N} \, \cos \left ( \frac{n (k+1/2) \pi}{N+1} \right ) \, \cos \left ( \frac{m (k+1/2) \pi}{N+1} \right ) \\&= \frac{1}{2}\Delta(m+n) + \frac{1}{2}\Delta(m-n),
\end{align}
where
\begin{equation}
    \label{eq:aliasing_delta}
   \Delta_N(\ell)  = 
\left\{
	\begin{array}{ll}
		1  & \text{if } \ell \, \mathrm{mod} \, (4(N+1)) = 0 \\
        -1  & \text{if } \ell \, \mathrm{mod} \, (4(N+1)) = 2(N+1) \\
		0 & \mbox{o.w. }
	\end{array}
\right.
\end{equation}
\end{lemma}

Substituting \Cref{eq:aliasing_delta} into \Cref{eq:quad_1D}:
\begin{align}
b_{n,N} &= \sum_{m=0}^M c_m\left(\frac{1}{2}\Delta_N(m+n) + \frac{1}{2}\Delta_N(m-n)\right) \\&= \frac{1}{N+1} \, \sum_{k=0}^N \cos(n\theta_k) \, p \left( \cos(m\theta_k) \right )  . \\    
\end{align}
This equation establishes a relationship between the aliased coefficients, the polynomial coefficients, and the polynomial values, which can be expressed as
\begin{equation}
    \vb_N = A \, \vc = F \, \vp,
\end{equation}
where $\vc$ is the vector of Chebyshev coefficients of the polynomial, $\vp$ contains polynomial values samples on a regular grid, $F$ corresponds to a size-$(N+1)$ Discrete Cosine Transform (DCT), and $A$ is sparse with entries given by
\begin{equation}
   A_{n,m} = \frac{1}{2}\Delta_N(m+n) + \frac{1}{2}\Delta_N(m-n).
   \label{eq:aliasing_matrix_1d}
\end{equation}

% Note that we obtained $b_{n,N}$ by computing the discrete cosine transform in frequency $n$ using $N+1$ equispaced points. Suppose these points are written as a vector $\vf$ then we have that:
% \begin{equation}
% F\vf = A\vc
% \end{equation}

% where $F$ is the linear operator applying different discrete cosine transforms on the samples to obtain the vector $b$.

\paragraph{Higher dimensions}
To see how this phenomenon impacts higher dimensions, we generalize the above relations from scalars in $[-1,1]$ to vectors $[-1,1]^D$. 
Chebyshev interpolation of a function $f:[-1,1]^D \to \mathbb{R}$ involves tensor-products of Chebyshev polynomials and replacing the single index $n$ with multi-index $\vn=(n_1, \hdots, n_D) \in \mathbb{N}^D$:
\begin{align*}
f(\vx) \approx p(\vx) = \sum_{\vn \in \mathcal{N}} c_\vn T_\vn(\vx), &\quad \vx \in [-1,1]^D, \;\\
&\quad c_\vn = \frac{2 - \delta_0(n)}{\pi}\int_{[-1,1]^D} \frac{T_\vn(\vx)f(\vx)}{\sqrt{\vector{1}-\vx^2}} d\vx,
\end{align*}
where $T_\vn(\vx) = T_{n_1}(x_1)T_{n_2}(x_2)\dots T_{n_d}(x_D)$, and $\mathcal{N}$ consists of the set of the multi-indices of the non-zero coefficients of the interpolant $p(\cdot)$. In general, for interpolation purposes, $\mathcal{N}$ corresponds to $ \mathcal{S}^1_{d,D}$ or $ \mathcal{S}^2_{d,D}$ (see \Cref{sec:preliminary}).

We apply the same change of variables in each dimension as in the one-dimensional case, producing
\begin{align}
p(\cos(\theta_1), ..., \cos(\theta_D)) ) &= \sum_{\vn \in \mathcal{N}} c_\vn\cos(\vn\vtheta), \quad \vtheta \in [0,\pi]^D, \\
c_\vn &= \int_{[0,\pi]^D} \cos(\vn\vtheta)f(\cos(\theta_1), ..., \cos(\theta_D)) d\vtheta
\label{eq:poly_approx_nd}
\end{align}
with $\vtheta = (\theta_1,\hdots, \theta_D) \in [0, \pi]^D$ and \
$\cos(\vn\vtheta) = \cos(n_1\theta_1)\cos(n_2\theta_2) \dots \cos(n_D\theta_D)$. Combining the two parts of \Cref{eq:poly_approx_nd}, we obtain
\begin{align}
c_\vn 
% & \sum_{\vm \in \mathcal{N} } c_\vm\int_{[0,\pi]^D} \cos(\vn\vtheta) \cos(\vm\vtheta) d\vtheta \\
%=& \sum_{\vm \in \mathcal{N}} c_\vm\int_{[0,\pi]^D} \left(\prod_{i=0}^D\cos(n_i\theta_i)\right)\left(\prod_{i=0}^D \cos(m_i\theta_i)\right) d\vtheta  \\
=& \sum_{\vm \in \mathcal{N} }c_\vm\prod_{i=0}^D\int_0^\pi \cos(n_i\theta_i)\cos(m_i\theta_i) d\theta_i.\label{eq:inner_product_nd} 
\end{align}

% MJM I think this isn't right...
%Discretizing the integral in \cref{eq:inner_product_nd} using $(N+1)^D$ equispaced points in $[0, \pi]^D$, we produce a tensor of samples $\vThetc_\vk$ indexed by the multi-index $\vk=(k_1, \hdots, k_D)$:
%\begin{equation}
%    \vThetc_\vk = \frac{\left(\vk + \frac{1}{2}\right)\pi}{\vN+1}, \quad \vk=\vzero,\hdots, \vN \in \{\vk\in \mathbb{N}^D\mid 0\leq k_i \leq N,\,\, i = 0, \hdots D\}.
%\end{equation}
Discretizing the integral in \Cref{eq:inner_product_nd}, using $N_i + 1$ points in each dimension $i$, provides the multi-dimensional analogue of \Cref{eq:quad_1D}:
\begin{align}
b_{\vn, \vN} &=
\sum_{\vm\in \mathcal{N}}
c_\vm
\prod_{i=0}^D
\left(\frac{1}{ N_i}
\sum_{k_i =0}^{N_i} 
\cos(n_i \, \theta_{i,k_i} ) \cos(m_i \, \theta_{i,k_i}) 
\right) \\
&= \sum_{k_1=0}^{N_1}  ... \sum_{k_D=0}^{N_D}  \, \left ( \frac{1}{\prod_{i=1}^D N_i} \,  \prod_{i=0}^D \cos(n_i\theta_{i,k_i}) \right ) \, p \left (\cos(\theta_{1,k_1}), \,  ..., \, \cos(\theta_{D,k_D}) 
\right ).
\label{eq:quad_nd}
\end{align}
In particular, note that if we sample on a grid where $N_i \geq \max_{n\in \mathcal{N}} (\max_i n_i)$ as is done traditionally, then $b_{\vn,\vN} = c_\vn$. However, as the number of samples increases exponentially in higher dimensions, this traditional approach becomes prohibitively expensive computationally. To alleviate this issue, we perform a form of under-sampling, leveraging the aliasing result of \Cref{lem:aliasing}; i.e., we use grids with $N_i$ points per dimension $i$, chosen in such a way that the total size $\prod_{i=1}^D N_i$ is roughly equal to $\lvert \mathcal{N} \rvert$ (see \Cref{sec:sampling_schedule} for details). Following \Cref{lem:aliasing}, this results in a linear system of the form
\begin{align}
b_{\vn, \vN} &= \sum_{\vm \in \mathcal{N}} 
c_\vm
\prod_{i=0}^D
\left(
\frac{1}{2}\Delta_{N_i}(m_i+n_i) + \frac{1}{2}\Delta_{N_i} (m_i-n_i)
\right)\\
&= \sum_{k_1=0}^{N_1}  ... \sum_{k_D=0}^{N_D}  \, \left (  \prod_{i=0}^D \cos(n_i\theta_{i,k_i}) \right ) \, p \left (\cos(\theta_{1,k_1}), \,  ..., \, \cos(\theta_{D,k_D}) 
\right ),
\label{eq:quad_nd_matrix} 
\end{align}
which may be written as
\begin{equation}
b_{\vN}  = A_{\mathcal{N}, \vN} \vc = F_{\vN} \vp_\vN.
\label{eq:multi_dim_dct}
\end{equation}
Here, $\vc$ is the vector of Chebyshev coefficients of the polynomial with multi-index corresponding to elements of $\mathcal{N}$, $\vp_\vN$ contains polynomial values from the non-uniform grid with $N_i+1$ samples per dimension, $F_{\vN}$  corresponds to a non-uniform Discrete Cosine Transform (DCT) with $N_i+1$ samples per dimension, and $A_{\mathcal{N}, \vN}$ is a size-$O(\lvert \mathcal{N} \rvert)  \times \lvert \mathcal{N} \rvert $ matrix with rows corresponding to elements of the non-uniform grid and columns corresponding to elements of the set $\mathcal{N}$ with entries
\begin{equation}
   A_{\vn, \vm}= \prod_{i=1}^D \, A^{(i)}_{n_i, m_i}  = \prod_{i=0}^D \left(\frac{1}{2}\Delta_{N_i}(m_i+n_i) + \frac{1}{2}\Delta_{N_i}(m_i-n_i)\right).   \label{eq:aliasing_matrix_nd_diff}
\end{equation}
In particular, the size of this system is \emph{linear} in the number of Chebyshev coefficients ($\lvert \mathcal{N} \rvert$), which leads to computational efficiency. Further, the matrix $A$ can be computed efficiently by considering the products of non-zero elements of lower-dimensional sparse aliasing matrices only.

\subsection{\texorpdfstring{$L$}{L}-grids}
With the system matrix on a single tensor product grids defined, we can apply the same construction to a collection of different grids. 
Let $\vNbar^{(\ell)} = (N_1^{(\ell)},\hdots, N_D^{(\ell)})$ for $\ell=1, \hdots, L$ be the $\ell^{th}$ choice of tensor-product discretization. 
The $L$ discretizations resulting from discretizing each dimension via  \Cref{eq:chebyshev_grid} for each $\ell$ is called an \textit{$L$-grid}, which is completely determined by the sampling rates $\vNbar^{(\ell)}$.

For each grid within a given $L$-grid, we compute the matrix $A^{(\ell)}$ using \Cref{eq:aliasing_matrix_nd_diff} and concatenate them into a final full matrix $A$:
\begin{equation}
A = \begin{bmatrix}
    A^{(1)} \\
    A^{(2)} \\
    \vdots \\
    A^{(L)} 
\end{bmatrix}
\label{eq:L-grid_A_mat}
\end{equation}

Combining \Cref{eq:L-grid_A_mat} with \Cref{eq:multi_dim_dct}, we obtain the expression in \Cref{eq:main_linear_system}. From this, the FCT computes the DCT on the left side of \Cref{eq:main_linear_system} using the function samples on the entire $L$-grid, followed by solving a least-squares problem for the unknown polynomial coefficients $\vc$.

\subsection{Sampling and Matrix Construction}
\label{sec:sampling_schedule}
With the system in \Cref{eq:main_linear_system}, we need to determine an $L$-grid, along with an appropriate value of $L$, such that \Cref{eq:main_linear_system} is well-conditioned and invertible. 
% A minimum requirement to guarantee the invertibility of $A$ is to choose an $L$-grid such that
% \begin{equation}
%     \sum_{\ell=1}^L \left ( \prod_{i=1}^D N_i^{(\ell)} \right ) \geq N.
% \end{equation}
% This is clearly too loose of a constraint: setting $N^{(\ell)}_i = 1$ will always insufficiently resolve $f$ in the $i$th dimension.

To ensure invertibility, we adopt a randomized approach by randomly assigning sampling rates along each dimension $i$ in such a way that
\begin{equation}
    \prod_{i=1}^D N_i^{(\ell)} \approx N, \quad \ell=1,\hdots L
\end{equation}
as detailed in \Cref{alg:sampling_rate_selection}.

With our randomly selected $\vNbar$, we sample $f$ on the associated $L$-grid and form the matrix $A$ in \Cref{eq:main_linear_system}. 
The goal for a sufficiently large value of $L$ is to sample enough points to produce a well-conditioned least-squares system to recover the Chebyshev coefficients. To determine $L$, we iteratively sample a new random grid, construct the corresponding matrix $A^{(\ell)}$, and append it to the matrix computed for the previous $\ell-1$ grids.
This process terminates when $A$ is full-rank with a condition number lower than some provided bound $\kappa$; we typically choose $\kappa< 10^4$.
This process is shown in \Cref{alg:construction_matrix_A}.
\begin{algorithm}
\caption{Sampling Rate Selection}
\label{alg:sampling_rate_selection}
\begin{algorithmic}[1]
\STATE{INPUT: $N$ (number of nonzero coefficients), $D$ (problem dimension), $d$ (degree)}
\STATE{OUTPUT: $(N_1, ..., N_D)$ (tuple of dimension-wise sampling rates)}
\STATE {Compute random permutation $(i_1, ..., i_D)$ of tuple $(1,...,D)$.}
\FOR{$j = 1,...,D$}
\STATE {Pick $N_{i_k}$ uniformly at random in $\{1, 2, ..., d+1\}$.}
\IF{$\prod_{k=1}^j N_{i_j} > N$}
\STATE {Set $N_{i_k} = 1$ for $k = j+1 ,...,D$.}
\STATE {break}
\ENDIF
\ENDFOR
\STATE {RETURN: $(N_1, ..., N_D)$}
\end{algorithmic}
\end{algorithm}

\begin{algorithm}
\caption{Construction of Aliasing Matrix $A$}
\label{alg:construction_matrix_A}
\begin{algorithmic}[1]
\STATE {INPUT: $\mathcal{N}$ (list of multi-indices of non-zero coefficients), $\kappa_0$ (maximum condition number)}
\STATE {OUTPUT: Matrix $A$ of \Cref{eq:main_linear_system} (full-rank), number of sampled grids $L$ , discretization pattern of $L$-grid $\mathfrak{N}$}
\STATE{Let $N = |\mathcal{N}|
$, $A = \emptyset$,$\mathfrak{N} = \emptyset$,  $L=0$.}
\WHILE{$\kappa(A) > \kappa_0$}
\STATE {$L \gets L + 1$}
\STATE {Compute $(N^{(L)}_1, ..., N^{(L)}_D)$ using \Cref{alg:sampling_rate_selection}}
\STATE {Compute $A^{(L)}$ using aliasing identity in \Cref{eq:aliasing_delta} and $\mathcal{N}$}
\STATE {$A \gets [A ; A^{(L)}]$}
\STATE {$\mathfrak{N} \gets [ \mathfrak{N}, (N^{(L)}_1, ..., N^{(L)}_D)]$}
\ENDWHILE
\STATE {RETURN: $A$, $L$, $\mathfrak{N}$}
\end{algorithmic}
\end{algorithm}

We have performed a variety of numerical experiments with the sampling scheme of \Cref{alg:sampling_rate_selection}. Empirically, choosing $L = c D$ for some $c \geq 2$ is sufficient to guarantee invertibility and well-conditioning of the matrix $A$ of \Cref{eq:main_linear_system}.
\paragraph{Precomputation of $L$-grids}
Since \Cref{alg:sampling_rate_selection} depends only on $\mathcal{N}$ and is independent of the sample values, we can precompute and store the sampling rates $\vNbar$ of an appropriate $L$-grid as well as the matrix $A$ for various values of $N$, $D$, $d$ and a given set of multi-indices $\mathcal{N}$ (e.g., $\mathcal{S}^1_{d,D}$ or $\mathcal{S}^2_{d,D}$).
% The matrix $A$ can also be stored, but since $A$ is fully determined by $\vNbar$, which is exactly $D$ integers and cheap to store. 
This makes \Cref{alg:construction_matrix_A} a precomputation with no runtime cost in practical applications.

\subsection{Full algorithm}
\label{sec:algorithm}
We now give a description of the full FCT algorithm.
Given a target function $f(\cdot)$, we first randomly sample an $L$-grid and construct the corresponding sparse matrix $A$ using  \Cref{alg:construction_matrix_A} and  \Cref{alg:sampling_rate_selection}.
Then, we sample the function $f(\cdot)$ at the discretization points on each grid in the chosen $L$-grid.
We denote the vector containing these samples on a single grid $\vf^{(\ell)}$ and the stacked vector of all values $\vf = [\vf^{(1)}; \hdots; \vf^{(L)} ]$.
We then form the linear system in \Cref{eq:main_linear_system} and solve for the Chebyshev coefficients $\vc$ via least-squares:
\begin{equation}
    \vc^* = \mathrm{argmin}_\vc ||(F \vf - A \vc)||_2^2.
\end{equation} 
We solve this least-squares problem with the CG method. Since we can efficiently compute Type-II DCTs with the Fast Fourier Transform (FFT) and can efficiently compute matrix-vector products $A\vc$ due to the sparsity of $A$, we can solve this problem rapidly.  Combined with the favorable conditioning of $A$, our approach allows the CG method to converge rapidly to a solution.
The scheme is summarized in \Cref{alg:sFCT}.

\begin{algorithm}
\caption{Sparse Fast Chebyshev Transform}
\label{alg:sFCT}
\begin{algorithmic}[1]
\STATE Input: target function $f$, degree $d$, dimension $D$, multi-index norm $s$
\STATE Output: Chebyshev coefficients $\vc^*$ approximating $f$
\STATE $\kappa_0 \gets 10^4$
\STATE $A, L, \mathfrak{N} \gets $ \Cref{alg:construction_matrix_A}($\mathcal{S}^s_{d,D}$, $\kappa_0$)
\FOR{$\ell = 1,\hdots,L$}
\STATE $\vf^{(\ell)} \gets$ Sample $f$ at $\mathfrak{N}_\ell$
\STATE $\mathcal{F}^{(\ell)} = F^{(\ell)} \, \vf^{(\ell)}$ (apply Type-II DCTs to $\vf^{(\ell)}$)
\ENDFOR
\STATE $\mathcal{F} = [\mathcal{F}^{(1)};\hdots; \mathcal{F}^{(L)}]$
\STATE Solve least-squares problem: $\vc^* = \mathrm{argmin}_\vc ||(\mathcal{F} - A \vc)||_2^2$ using CG.
\RETURN $\vc^*$
\end{algorithmic}
\end{algorithm}

\subsection{Scaling \& Complexity}
\label{sec:complexity}
\paragraph{Precomputation complexity}
The primary work of the precomputation phase of the FCT involves forming the matrix $A$ using \Cref{alg:construction_matrix_A}. 
Computing the sampling rate for a given block of $A$, which calls \Cref{alg:sampling_rate_selection}, requires $O(D)$ time to complete.
For a given sampling rate $(N_1^{(\ell)}, N_2^{(\ell)}, \hdots, N_D^{(\ell)})$, we then compute the non-zero elements of $A^{(\ell)}$ by using the nonzeros of one dimensional aliasing matrices as in \Cref{eq:aliasing_matrix_nd_diff};  since each one dimensional aliasing matrix has $O(1)$ zeros per dimension per column, computing all of the zeros results in $O(ND)$ computational complexity.
Further, because the loop in \Cref{alg:construction_matrix_A} repeats $L$ times, we can construct $A$ in $O(NDL)$ time (the construction of the $\ell^{th}$ aliasing matrix dominates the cost of each loop iteration. It is also worth noting that, as the matrix $A$ is sparse, the memory complexity is $O(NL)$.

 % Since the columns of $A^{(\ell)}$ are Kr\"onecker product of $D$ one-dimensional sparse aliasing matrices each with $O(N_i^{(\ell)})$ non-zero entries (\Cref{eq:aliasing_matrix_nd_diff}), we need to compute the product of $D$ elements for non-zero entries only, which can be done in $O(D\prod_{i=1}^D N_i^{(\ell)}) = O(ND)$ time.

\paragraph{FCT complexity}
There are two main steps in \Cref{alg:sFCT}: (1) the Type-II DCTs applied to the function samples and (2) solving the least-squares problem with the CG method.
The adjoint Type-II DCTs applied on on $L$ different $D$-dimensional grids (containing $O(N)$) points can be computed in $O(NL\log(N))$  time using FFTs \cite{ahmed1974discrete}.
The least-squares problem can be solved with the CG method in $n_{CG}$ iterations with each iteration requiring a sparse matrix-vector multiply in $O(NL)$ time, resulting in an overall computational complexity of $O(n_{CG}NL)$.
Note that if the condition number of $A$ is not too large in value, resulting in a {\em well-conditioned} system, $n_{CG}$ is small \cite[Chapter 10.2]{golub2013matrix}, and the CG method converges rapidly.
 We discuss this further in the following section, but in the case of a well-conditioned $A$,the complexity of FCT is $O(NL\log(N))$.

In summary, the precomputation phase computational complexity is $O(NDL)$, while applying FCT to $f$ results in $O(NL\log(N))$ complexity.

\section{Numerical Results}
\label{sec:numerical_results}
In this section, we show results in terms of scalability and accuracy for FCT.  We begin with a brief background introduction on numerical approximation and error in \Cref{sec:approx_and_details}, followed by a discussion on the methods used for comparison in \Cref{sec:numerical_comparisons} as well as an introduction to the different types of functions studied for scaling and accuracy in \Cref{sec:numerical_families}.  We briefly note some details on the numerical implementation in \Cref{sec:implementation}.  We complete our numerical discussions by showing the specific experiment results on scaling and accuracy in \Cref{sec:numerical_experiment_results} for FCT against the methods in \Cref{sec:numerical_comparisons} on the functions in \Cref{sec:numerical_families}.  
%Finally, we provide a discussion for our empirical observations on choosing grid-sampling parameters in \Cref{sec:empirical_L}.

\subsection{Numerical Approximation and Error Discussion} 
\label{sec:approx_and_details}
The error introduced from the polynomial approximation of a smooth function $f$ evaluated at Chebyshev points is a well-studied problem \cite{trefethen2019approximation}. 
On a tensor product Chebyshev grid in high dimensions, if $f$ is analytic in an Bernstein $\rho$-hyperellipse containing $[-1,1]^D$, then an interpolating polynomial of Euclidean degree $k$ will converge geometrically with rate $O(\rho^{-k})$ \cite{trefethen2017multivariate}.
In the context of FCT, we are solving a least-squares problem for an approximating polynomial to $f$ on an $L$-grid instead of interpolating $f$ on a single $D$-dimensional grid.
The quadrature rule used in \Cref{eq:quad_nd} to compute $c_\vn$ in \Cref{eq:inner_product_nd} is exact due to the orthogonality of cosines; however, when $f$ is underresolved in a certain dimension as a result of the randomized grid sub sampling, greater care must be taken.
Suppose $f$ is a polynomial with bounded total degree $k$, but whose leading order term has a multi-index $\vk = (k_1, k_2, \hdots, k_D)$ with $\sum_{i=1}^D k_i \leq N$. We have observed empirically that as long as $L$ scales on the order $O(D+d)$, where $D$ is the dimension of the ambient space and $d$ is the maximum degree we wish to interpolate, the resulting aliasing matrix $A$ has full column rank with high probability; indeed, experiments show this scaling for $L$ holds for all numerical results considered below in the following sections. In this case, we can obtain an an accurate Chebyshev expansion approximating $f$ using only a constrained set of polynomials in our basis from $\mathcal{N}$.

With this background going forward, we next discuss the methods we employ for evaluating the accuracy and scaling benefits of FCT in \Cref{sec:numerical_comparisons}.

% Specifically, we use the FCT to compute the coefficients of bounded total degree polynomial (i.e. in $\PDdone$); in \Cref{sec:results_sparse}, we apply FCT to a modified sparse polynomial recovery problem. 
% FCT exhibits quasi-linear scaling in each case and outperforms both competing approaches.
% In \Cref{sec:accuracy_experiments}, we demonstrate the accuracy and performance of FCT on a concrete interpolation problem. 

\subsection{Methods for Numerical Comparison}
\label{sec:numerical_comparisons}
We compare FCT with two standard approaches for performing Chebyshev interpolation: (1) the baseline DCT-based interpolation approach and (2) a least-squares approach based on a randomized interpolation matrix, which we call Randomized Least-Squares Interpolation (RLSI) with the following specific details:
\begin{itemize}
\item {\bf DCT approach}: we sample the target function $f$ on a $D$-dimensional tensor product Chebyshev grid and apply the DCT to compute the Chebyshev coefficients.
This method has a complexity of $O(Dd^D\log d)$.
\item {\bf RLSI approach}: consider the general least squares approach as introduced in \Cref{sec:preliminary}.  For a Randomized Least Squares (RLSI) approach, we construct the matrix $B$ in \Cref{eq:least_squares_matrix} to solve \Cref{eq:least_squares_opt} by uniformly sampling $CN$ points from a full uniform tensor-product grid of size $d^D$ (without explicitly forming it). We further choose $C=1.2$ to ensure that $B$ is invertible with bounded condition number less than $\kappa_{\mathrm{RLSI}} = 10^4$; this mirrors the conditions imposed in \Cref{alg:construction_matrix_A} for a fair comparison. 
The resulting least-squares problem is solved numerically using CG, available through the \texttt{Eigen} \cite{eigenweb}. This method has complexity $O(\vert \mathcal{N} \rvert^2)$.
\end{itemize}
Note that for consistency, we use a tolerance of $10^{-3}$ for both RLSI and FCT.
We choose $L=3D$ for FCT for all experiments unless otherwise specified.

%We first compare these algorithms on several function approximation problems. 
%Then, we will compare FCT with the DCT-based approach and RLSI to a structured sparse coefficient recovery problem, where non-zero polynomial indices are known in advance but not their values.

To demonstrate the scalability and efficiency of the FCT algorithm versus the DCT-based and RLSI-based approaches, we compare against function samples for three different challenging families of functions as detailed in \Cref{sec:numerical_families}.

\subsection{Numerical Experiments Function Families}
\label{sec:numerical_families}
For our numerical experiments, we used three families of functions as detailed below; these three function families are studied extensively in \Cref{sec:numerical_experiment_results} against the methods as detailed above in \Cref{sec:numerical_comparisons}.  These function families are chosen for their challenging complexity in approximation, interpolation, and scalability by FCT, DCT-based, and RLSI methods, providing a clear and fair comparison of all three methods.

\subsubsection{Function Class \#1 (Runge Function)}
\label{sec:numerical_runge}
The first family $f_1(\cdot)$ can be parameterized as:
\begin{equation*}
    f_1(\vx) = \frac{1}{1 + 10\|\vx\|_2^2}.
\end{equation*}

This function is challenging to approximate accurately due to the nearby pole in the complex extension of $f_1$ from $\mathbb{R}^D$ to $\mathbb{C}^D$ located at $\pm \frac{i}{\sqrt{10}}\vector{e}_i$ for $i^{th}$ unit vector, $\vector{e}_i$.  Results for experiments using this function can be seen in \Cref{sec:numerical_runge_res}.
\subsubsection{Function Class \#2 (Oscillatory Function)}
\label{sec:numerical_oscillatory}
The second family of functions, $f_2(\cdot)$, are parameterized as:
\begin{equation}
    f_2(\vx) = \sin(3\cos(3\exp(\|\vx\|_2^2))) + \exp(\sin(3(\vx \cdot \vector{1}))).
\end{equation}

These functions are analytic, but they exhibit a complex behavior in the region $[-1,1]^D$. In this case, accurate interpolation requires high degree, and therefore many coefficients, making this family challenging computationally. Results for experiments using this function can be seen in \Cref{sec:numerical_osc_res}.

%as shown in Figure \ref{fig:f2_2d}. 

% \begin{figure}[!htp]
% \captionsetup[subfigure]{labelformat=empty}
%      \centering
%      \begin{subfigure}[b]{0.495\textwidth}
%          \centering
%          \includegraphics[width=\textwidth]{accuracy_function.pdf}
%      \end{subfigure}
%         \caption{\small Two-dimensional rendering of the function $f_2(\cdot)$ in the region $[-1,1]$. }
%         \label{fig:f2_2d}
% \end{figure}
\subsubsection{Function Class \#3 (Sparse Support Function)}
\label{sec:numerical_sparse_support}
Finally, we consider a third family of functions $f_3(\cdot)$, which have \emph{sparse support}; i.e., $\mathcal{N}$ has little structure and is a sparse subset of $S^\infty_{d,D}$. To generate elements of this family, we fix a dimension $D$ and degree $d$ and randomly select $N$ multi-indices within $S^\infty_{d,D}$. We then assign random values in $[-1,1]$ to the coefficients corresponding to the chosen multi-indices.  Results for experiments using this function can be seen in \Cref{sec:numerical_sparse_res}.

\subsection{Numerical Implementation}
\label{sec:implementation}
We have implemented the FCT algorithm described in \Cref{sec:algorithm} in \texttt{C++}.
We use the \texttt{fftw} package \cite{frigo1998fftw} for its optimized DCT implementation and \texttt{Eigen} \cite{eigenweb} for its linear algebra routines, including its CG method.
 All experiments were performed with a single core on a machine with an AMD EPYC 7452 processor and 256GB of RAM.

\subsection{Numerical Experiment Studies on Scaling and Accuracy}
\label{sec:numerical_experiment_results}
Using our FCT implementation, we provide the specific results of our numerical experiments against the methods detailed in \Cref{sec:numerical_comparisons} (i.e., DCT-based and RLSI) on the three function families described in \Cref{sec:numerical_families}.  We begin with Function Family \#1 (Runge Function $f_1(\cdot)$) as described in \Cref{sec:numerical_runge}.

\subsubsection{Runge Function (\texorpdfstring{$f_1(\cdot)$}{f1})}
\label{sec:numerical_runge_res}
We first focus on Function Family \#1 as described in \Cref{sec:numerical_runge} in order to study scaling of FCT against DCT-based and RLSI methods.

\paragraph{Scaling Comparison with $f_1(\cdot)$}
To compare the scaling behaviors of each method, we continue to focus on $f_1(\cdot)$ and construct polynomials in $\PDdone$ with random coefficients $c_\vn$ in $[-1,1]^D$ and attempt to recover them with FCT, the DCT-based approach, and RLSI. 

Each of these algorithms should recover the polynomials accurately with varying complexity.
We test polynomials of degree $d=3,6$ and in dimensions $D= 2, \hdots, 25$ with $N=|\mathcal{S}^1_{d,D}| = \binom{D+d}{d}$ coefficients.

\begin{figure}[!htp]
\captionsetup[subfigure]{labelformat=empty}
     \centering
     \begin{subfigure}[!htb]{0.495\textwidth}
         \centering
         \includegraphics[width=\textwidth]{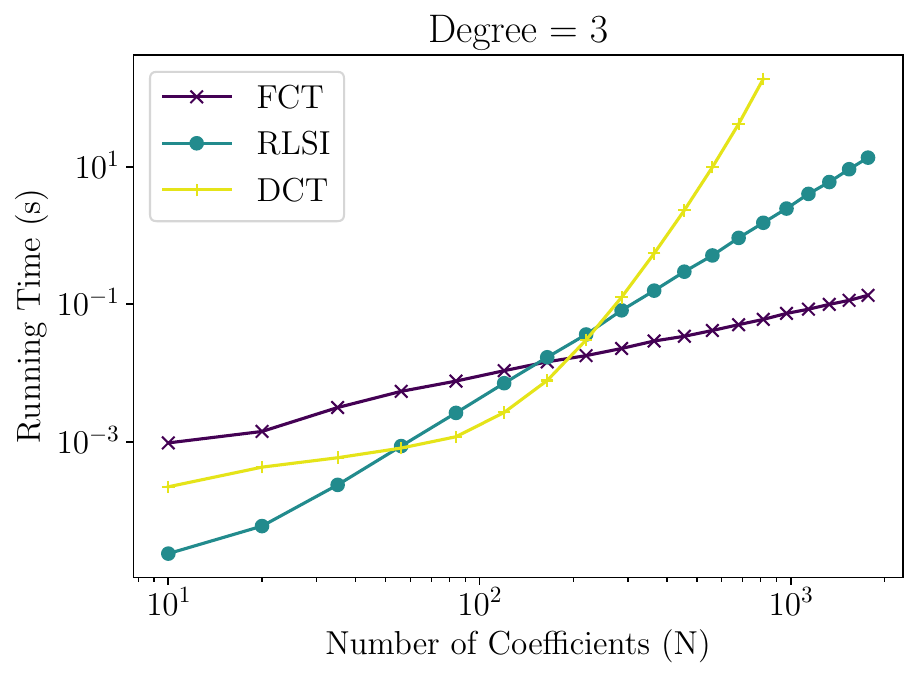}
         \label{fig:scaling_deg_3}
     \end{subfigure}
     \hfill
     \begin{subfigure}[!htb]{0.495\textwidth}
         \centering
         \includegraphics[width=\textwidth]{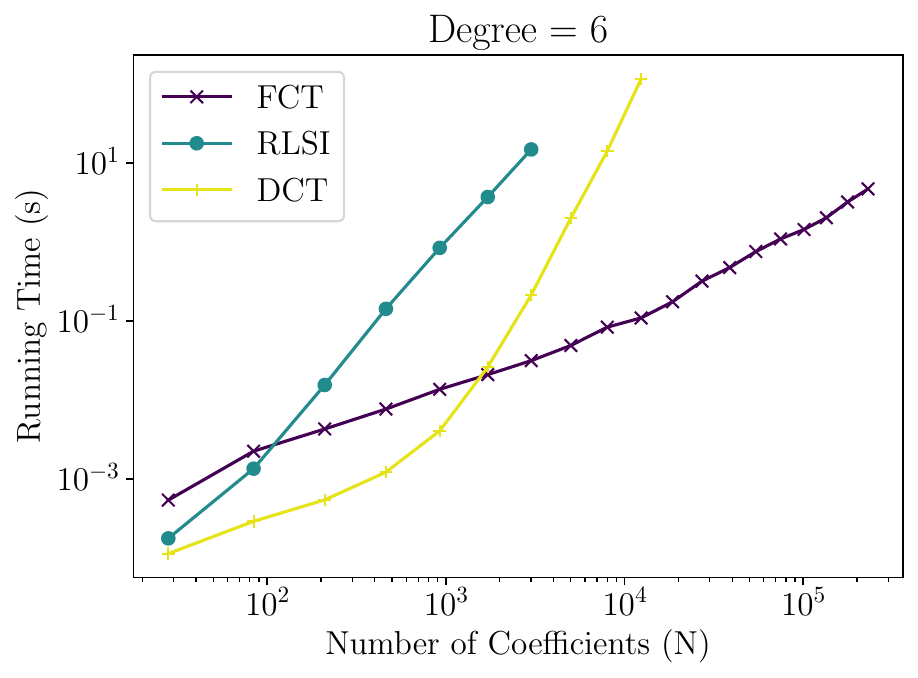}
         \label{fig:scaling_deg_6}
     \end{subfigure}
        \caption{\small Wall time versus number of coefficients $N=|\mathcal{S}^1_{d,D}| = \binom{D+d}{d}$ for the interpolation of polynomials with degree $d=3$ (left) and at $d = 6$ (right), and varying dimension. The FCT significantly outperforms the DCT and RLSI with a crossover size of roughly $N=2 \times 10^2$. Each point on each curve corresponds to dimension $D=2,\hdots,25$. 
        Curves with fewer points indicate the algorithm in question ran out of memory.}
        \label{fig:result_dense_support}
\end{figure}

%We summarize the results in \Cref{fig:result_dense_support} and \Cref{fig:result_dense_support_error}.  
In \Cref{fig:result_dense_support}, we report the wall time of each method in terms of the number of nonzero coefficients $N$ for a given value of $d$ and $D$.
We do not report timing results when the algorithm runs out of system memory.  We see that the runtime of the DCT approach grows exponentially with $D$, quickly running out of system memory even for $d=3$.
The RLSI approach fares better, exhibiting roughly $O(N^2)$ scaling as expected; however, it still runs out of memory when $d=6$ while FCT is able to recover the most polynomials in high dimensions.
It is worth highlighting that, in low dimensions, RLSI and the DCT approach are in fact faster than FCT.

\begin{figure}[!htp]
\captionsetup[subfigure]{labelformat=empty}
     \centering
     \begin{subfigure}[!htb]{0.495\textwidth}
         \centering
         \includegraphics[width=\textwidth]{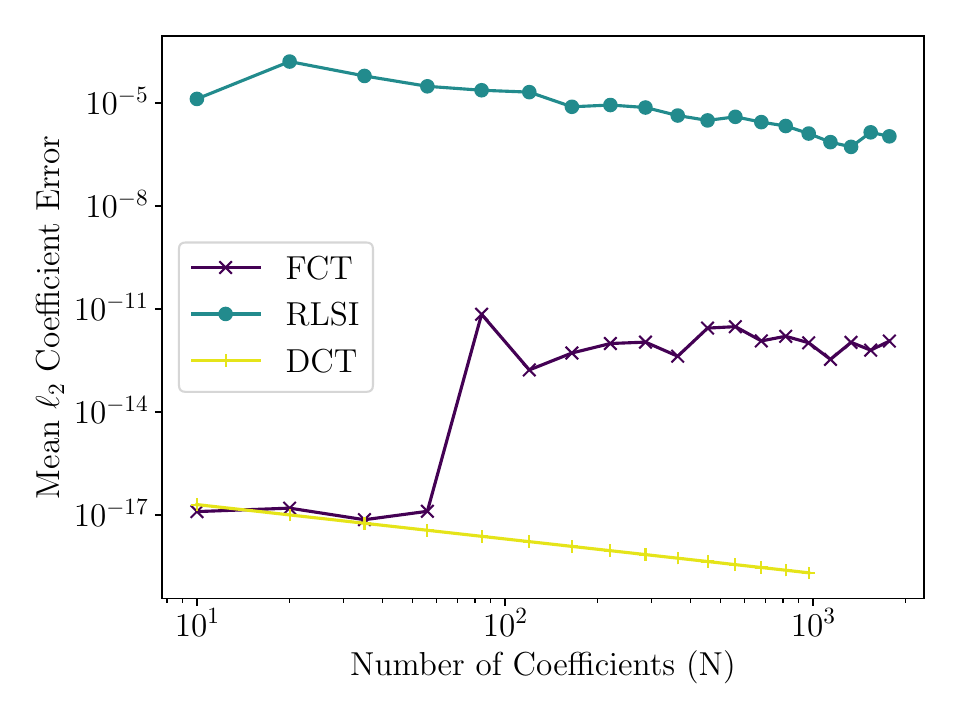}
         \label{fig:scalingerror_deg_3}
     \end{subfigure}
     \hfill
     \begin{subfigure}[!htb]{0.495\textwidth}
         \centering
         \includegraphics[width=\textwidth]{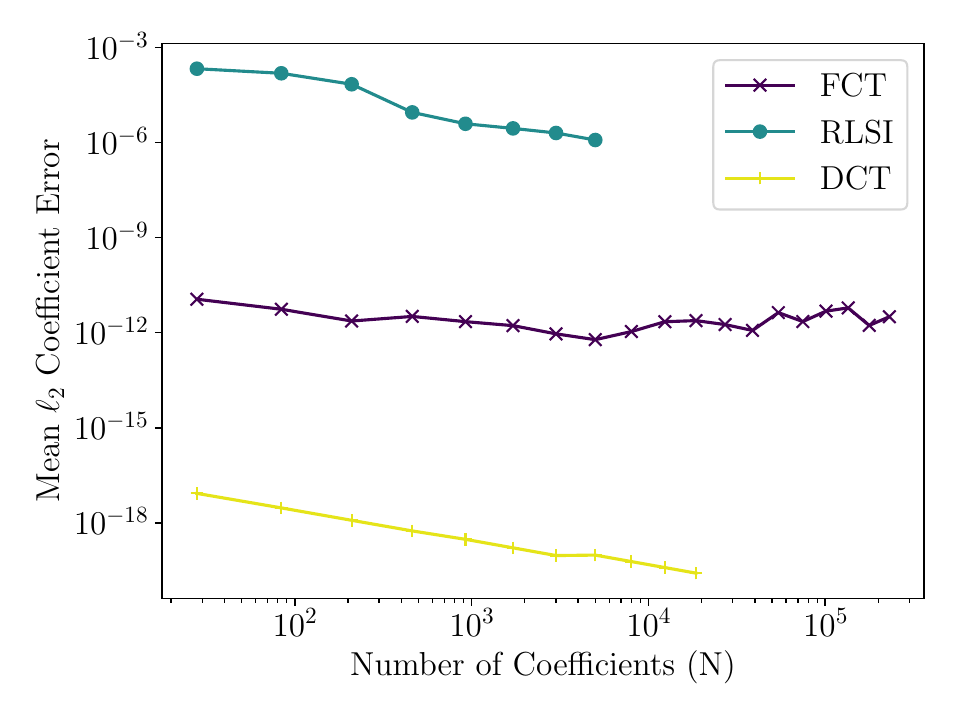}
         \label{fig:scalingerror_deg_6}
     \end{subfigure}
        \caption{\small Mean $\ell_2$ coefficient error versus number of coefficients $N=|\mathcal{S}^1_{d,D}| = \binom{D+d}{d}$ for the interpolation of polynomials with degree $d=3$ (left) and at $d = 6$ (right), and varying dimension. Each point on each curve corresponds to dimension $D=2,\hdots,25$.  }
        \label{fig:result_dense_support_error}
\end{figure}

\paragraph{Accuracy Comparison}  
In \Cref{fig:result_dense_support_error}, we compare the approximation accuracy of the three methods by providing the mean $\ell_2$ coefficient error: for a computed solution $\vc$ and known solution $\tilde{\vc}$, we report $\|\vc - \tilde{\vc}\|_2/N$.
Tests that run out of system memory in this range are not plotted. The target function for interpolation in this case was generated by randomly selecting Chebyshev coefficients between $[-1, 1]$ in the specified dimensions and degrees. 

The error of each method is essentially flat (negative slopes are due to increasing $N$).
The DCT-based approach exactly recovers the polynomial coefficients as expected, whereas 
the accuracy of RLSI is bounded by the CG error tolerance, since the system matrix is full-rank.
% Since the aliasing matrix of FCT has a rapidly decaying eigenspectrum, it is able to consistently recover a solution that is more accurate than the specified CG tolerance regardless of dimension and degree.
The jump in error for FCT in the degree case is caused by a phase change as we move from low condition number aliasing matrices for small problem instances, to ones with more moderate condition numbers.  We observe that FCT becomes more efficient than the DCT-based approach and RLSI at $N\approx 2\times 10^2$ for $d=3$ and $N\approx 2\times 10^3$ for $d=6$, i.e., $D=10$ for $d=3$ and $D=7$ for $d=6$.
% This is consistent with \Cref{table:sampling}, which shows that the DCT approach is faster in lower dimensions.

\subsubsection{Oscillatory Function  (\texorpdfstring{$f_2(\cdot)$}{f2})}
\label{sec:numerical_osc_res}
\begin{figure}[!htp]
\captionsetup[subfigure]{labelformat=empty}
     \centering
     \begin{subfigure}[b]{0.495\textwidth}
         \centering
         \includegraphics[width=\textwidth]{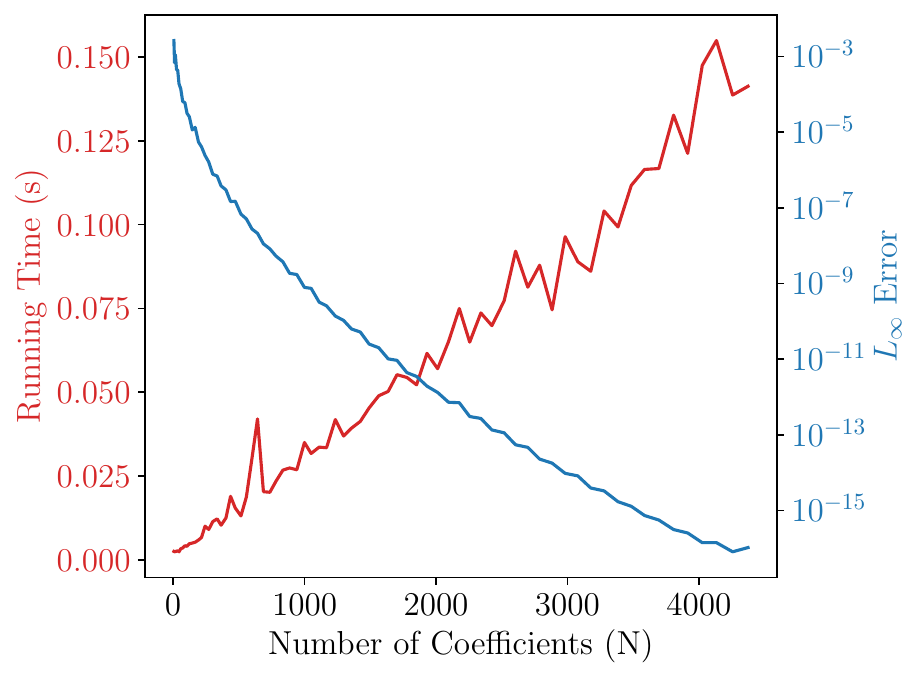}
     \end{subfigure}
     \hfill
     \begin{subfigure}[b]{0.495\textwidth}
         \centering
         \includegraphics[width=\textwidth]{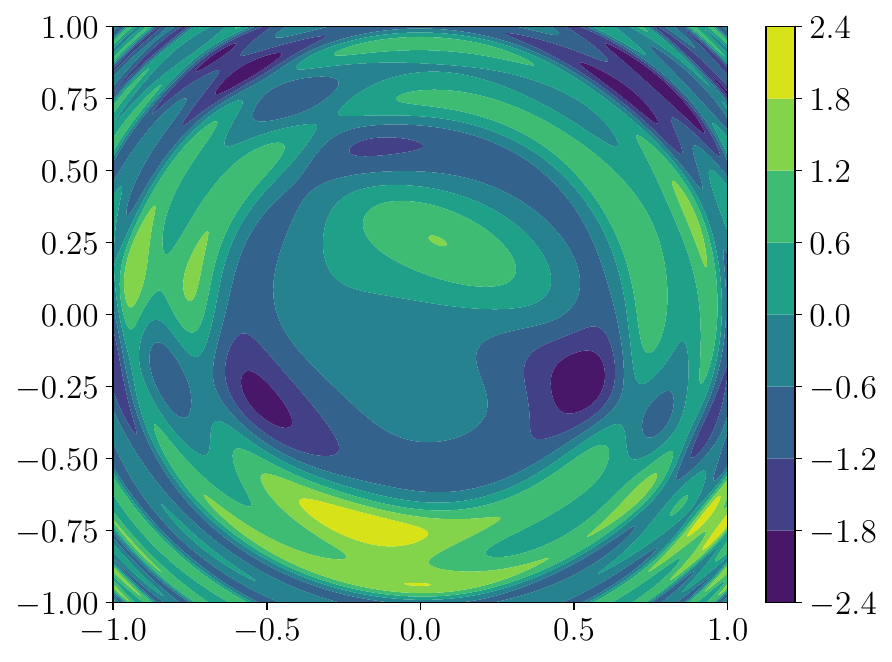}
     \end{subfigure}
        \caption{\small {\em Left}: FCT wall time and $L^\infty$ absolute error for $f_2$, with bounded Euclidean degree $d_E=2,\hdots, 50$; {\em Right}: a horizontal slice through the origin of $f_2$. We observe rapid error decay with increasing Euclidean degree, while the computational cost scales quasi-linearly.}
        \label{fig:accuracy}
\end{figure}

To investigate the accuracy of function approximations produced with FCT, we consider the function $f_2(\cdot)$ as described in \Cref{sec:numerical_oscillatory}, where we choose $D=5$.
A 2D slice of $f_2$ through the origin can be seen in \Cref{fig:accuracy}-Right.
We approximate $f_2$ with a polynomial of bounded Euclidean degree $d_E$; i.e.,  $N = |\mathcal{S}_{d,D}^2|$ and compare the $L^\infty$ absolute error.\footnote{To approximate the $L^\infty$ error, we compute the maximum error over a random sampling of 5000 points in $[-1,1]^D$.}.
We choose a CG tolerance at machine precision to ensure that all error is solely due to FCT.

In \Cref{fig:accuracy}-Left, we plot the relative error and wall time of FCT in terms of the number of coefficients $N$, determined by a bounded Euclidean degree.
We see that the $L^\infty$ error decays rapidly with $N$, and therefore with $d_E$, as expected from \cite{trefethen2017multivariate}.
By $d_E=50$ ($N\approx 5000$), we have resolved $f_2$ to machine precision, which is a behavior consistent with a standard Chebyshev interpolation scheme.
We also observe quasi-linear scaling in wall time with increasing $N$, requiring only a $1/10 s$ to recover 5000 coefficients to machine precision.

%and apply the FCT and RLSI algorithms to solve for the unknown coefficients. 
%we randomly select $N$ multi-indices in dimension $D=100$ up to degree $d=50$. 

\subsubsection{Sparse Coefficient Recovery (\texorpdfstring{$f_3(\cdot)$}{f3})}
\label{sec:numerical_sparse_res}
We can use the FCT algorithm to approximate a function $f$ using an expansion consisting of specific subset of Chebyshev coefficients as described in \Cref{sec:numerical_sparse_support}. Provided we know the location of the nonzero coefficients in the expansion of our function of interest, FCT will recover these rapidly and with high accuracy.
Suppose we know that the target function $f$ is a polynomial with $N$ non-zero coefficients whose locations we know in advance.
In other words, we know the multi-indices $\vn^{(i_1)},\vn^{(i_2)},\hdots,\vn^{(i_N)}$ which have non-zero coefficients $c_{\vn^{(i_j)}}$.

\begin{figure}[!htbp] 
     \centering
\includegraphics[width=0.5\textwidth]{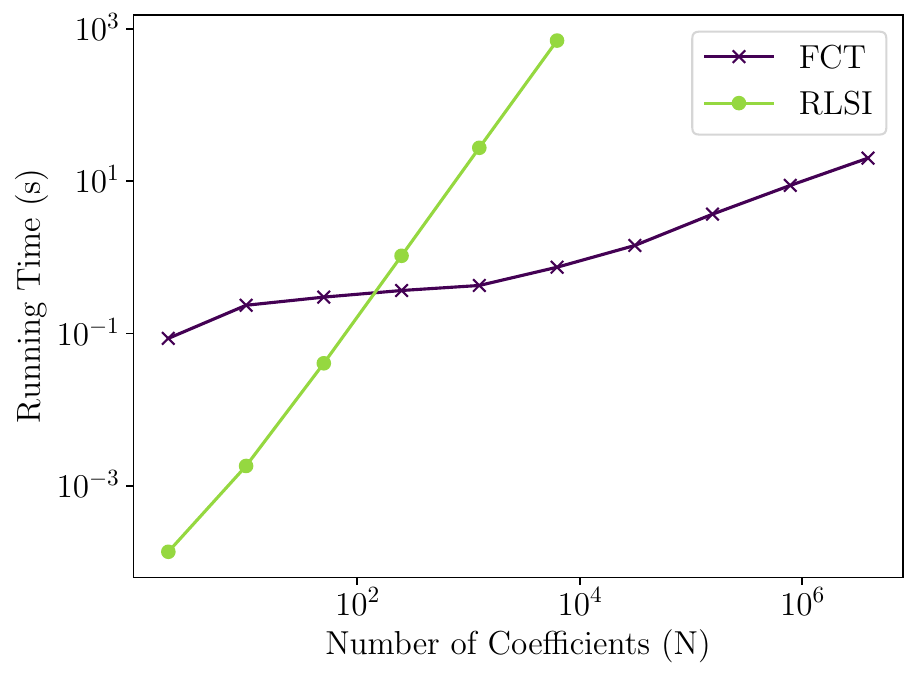}
    \caption{\small Comparison of the computation time of the FCT and RLSI algorithms for the recovery of the Chebyshev coefficients of a polynomials in dimension $D=100$ having an increasing number of sparse coefficients. 
    }
  \label{fig:wall_time_big_scaling}
\end{figure}

We compare the walltime for FCT and RLSI in \Cref{fig:wall_time_big_scaling} in dimension $D=100$. FCT recovers $N=10^7$ Chebyshev coefficients in $\sim 10$ seconds, while RLSI runs out of memory at $N=10^4$. Note that $A$ is constructed via equation \Cref{eq:aliasing_matrix_nd_diff}, and only contains columns corresponding to the nonzero multi-indices (i.e., only $N$ columns).

\section{Conclusions}
\label{sec:conclusions}
    In this paper, we have introduced the Fast Chebyshev Transform (FCT), an algorithm for computing the Chebyshev coefficients of a polynomial from the knowledge of the location of its nonzero Chebyshev coefficients with applications to Chebyshev interpolation.
    We have theoretically and empirically demonstrated quasi-linear scaling and competitive numerical accuracy on high-dimensional problems. 
    We have compared the performance of our method with that of commonly-used alternative techniques based on least-squares (RLSI) and the Discrete Cosine Transform (DCT).  We have further demonstrated significant speedups. 

    Moving forward, our most immediate next step is to improve the sampling scheme of FCT to further reduce the number of required function samples.
    We also plan to apply FCT to problems in polynomial optimization and apply the sampling scheme to more general transforms applied to smooth functions, such as the NUFFT (Non-Uniform Fast Fourier Transform) and Gaussian process regression \cite{greengard2022equispaced,greengard2004accelerating}.

\section{Data Availability Statement}
The data produced for the redaction of this paper will be made available under reasonable requests satisfying the policy of Qualcomm Technologies Inc. for the public release of internal research material.

\appendix
\section{Proofs} 
\subsection{Proof of aliasing identity}
\label{sec:lemma_proof}
% \lipsum[71]

\begin{lemma}
 Let $ 1 < N \in \mathbb{N}$. Then, 
 \begin{align}
    &\frac{1}{N} \, \sum_{k=0}^{N-1} \, \cos \left ( \frac{n (k+1/2) \pi}{N} \right ) \, \cos \left ( \frac{m (k+1/2) \pi}{N} \right )  &=\frac{1}{2} \Delta(m+n) + \frac{1}{2}\Delta(m-n) 
    % \frac{\cos( \pi (m+n) /2N) \,}{2} \Delta(m+n) + \frac{\cos( \pi (m-n) /2N) \,}{2}\Delta(m-n) 
\end{align}
where,
\begin{equation}
    \label{eq:aliasing_delta_app}
   \Delta(l)  = 
\left\{
	\begin{array}{ll}
		1  & \mbox{if } l \, \mathrm{mod} \, (4N) \equiv 0 \\
        -1  & \mbox{if } l \, \mathrm{mod} \, (4N) \equiv 2N \\
		0 & \mbox{o.w. }
	\end{array}
\right.
\end{equation}
\end{lemma}
\begin{proof}
Consider,
\begin{align}
    &\frac{1}{N} \, \sum_{n=0}^{N-1} \, \cos \left ( \frac{m (n+1/2) \pi}{N} \right ) \, \cos \left ( \frac{k (n+1/2) \pi}{N} \right ) \\
    &=\frac{1}{N} \, \sum_{n=0}^{N-1} \, \left [ \frac{1}{2} \cos \left ( \frac{(n+1/2) (m+k) \pi}{N} \right )  + \frac{1}{2}  \cos \left ( \frac{ (n+1/2) (m-k) \pi}{N} \right )  \right ] 
\end{align}
Then, consider,
\begin{align}
     \sum_{n=0}^{N-1} \, \cos \left ( \frac{(n+1/2) l \pi}{N} \right ) &= \frac{e^{i \frac{\pi}{2N} l}}{2} \, \sum_{n=0}^{N-1}e^{i \frac{\pi}{N} n \, l} + \frac{e^{-i \frac{\pi}{2N} l}}{2} \, \sum_{n=0}^{N-1}e^{-i \frac{\pi}{N} n \, l}  \\
    &= \left ( \frac{e^{i \frac{\pi}{2} l}}{2} + \frac{e^{-i \frac{\pi}{2} l}}{2}  \right ) \, \frac{\sin(\pi l/2) }{\sin( \pi l /(2N) )}
    \\
    &= \cos( \pi l /2) \, \frac{\sin(\pi l/2) }{\sin( \pi l /(2N) )}
    \\
    &=  \frac{1}{2} \, \frac{\sin(\pi l) }{\sin( \pi l /(2N) )}
\end{align}
since (Dirichlet kernel),
\begin{equation}
    \sum_{n=0}^{N-1} e^{i n \, x}  = e^{i \frac{N-1}{2} x} \, \frac{\sin(N x/2) }{\sin( x/2)}
\end{equation}
and since,
\begin{equation}
    \cos( \pi l /2) \, \sin(\pi l/2) = \frac{1}{2} \sin(\pi l) + \frac{1}{2} \sin(0)= \frac{1}{2} \sin(\pi l) 
\end{equation}
Therefore,
\begin{align}
    &\frac{1}{N} \, \sum_{n=0}^{N-1} \, \left [ \frac{1}{2} \cos \left ( \frac{(n+1/2) (m+k) \pi}{N} \right )  + \frac{1}{2}  \cos \left ( \frac{ (n+1/2) (m-k) \pi}{N} \right )  \right ] \\
    &= \frac{1}{2}  \, \left ( \frac{1}{2N}\, \frac{\sin(\pi (m+k)) }{\sin( \pi (m+k) /(2N) )} + \frac{1}{2N}\, \frac{\sin(\pi (m-k)) }{\sin( \pi (m-k) /(2N) )} \right )  \\
    &= \frac{1}{2} \, \Delta(m+k) + \frac{1}{2} \, \Delta(m-k) 
\end{align}
and the result follows after defining,
\begin{equation}
  \Delta(l) :=   \frac{1}{2N}\, \frac{\sin(\pi l) }{\sin( \pi l /(2N) )} = 
\left\{
	\begin{array}{ll}
		1  & \mbox{if } l \, \mathrm{mod} \, (4N) \equiv 0 \\
        -1  & \mbox{if } l \, \mathrm{mod} \, (4N) \equiv 2N \\
		0 & \mbox{o.w. }
	\end{array}
	\right.
\end{equation}
for all $l \in \mathbb{Z}$.
\end{proof}

% \section*{Acknowledgments}

\bibliographystyle{siamplain}
\bibliography{references}
\end{document}

%% file: shared.tex
\usepackage{lipsum}
\usepackage{amsfonts}
\usepackage{graphicx}
\usepackage{epstopdf}
\usepackage{algorithmic}
\ifpdf
  \DeclareGraphicsExtensions{.eps,.pdf,.png,.jpg}
\else
  \DeclareGraphicsExtensions{.eps}
\fi

\newsiamremark{remark}{Remark}
\newsiamremark{hypothesis}{Hypothesis}
\crefname{hypothesis}{Hypothesis}{Hypotheses}
\newsiamthm{claim}{Claim}

\title{A Sparse Fast Chebyshev Transform for High-Dimensional Approximation}

\author{Dalton Jones\thanks{Qualcomm AI Research. Qualcomm AI Research is an initiative of Qualcomm Technologies, Inc.} \and Pierre-David Letourneau\footnotemark[2] \and Matthew J. Morse\footnotemark[2] \and M. Harper Langston\footnotemark[2]}
\usepackage{amsopn}